\documentclass[11pt]{amsart}

\newtheorem{thm}{Theorem}
\newtheorem{lemma}[thm]{Lemma}
\newtheorem{prop}[thm]{Proposition}
\theoremstyle{definition}
\newtheorem{defn}{Definition}
\newtheorem{remark}[thm]{Remark}

\newcommand{\sinsq}{\sin^2 \theta}
\newcommand{\np}{{\bf n}}
\newcommand{\bt}{{\tilde{B}}}
\newcommand{\yt}{\tilde{y}}
\newcommand{\del}{\partial}
\newcommand{\law}{{\mathcal L}}

\newcommand{\what}{\widehat}
\newcommand{\Y}{{\bf Y}}
\newcommand{\R}{{\bf R}}
\newcommand{\C}{{\bf C}}
\newcommand{\Z}{{\bf Z}}
\newcommand{\half}{\frac{1}{2}}

\newcommand{\ignore}[1]{}

\begin{document}
\title[The Drunkard's Walk on the Sphere]{Discrepancy convergence for
the drunkard's walk on the sphere}

\author{Francis Edward Su}

\dedicatory{Manuscript version, February 2001}

\address{Department of Mathematics\\ Harvey Mudd College\\
Claremont, CA  91711}

\email{su@math.hmc.edu \qquad {\it web:} http://www.math.hmc.edu/$^\sim$su/}

\thanks{
The author thanks Persi Diaconis for suggesting this problem, Ken
Ross for useful feedback, and the Cornell School of Operations
Research for hospitality during a sabbatical where this was
completed.}

\subjclass{Primary 60B15; Secondary 43A85}

\keywords{discrepancy, random walk, Gelfand
pairs, homogeneous spaces, Legendre polynomials}

\begin{abstract}
We analyze the drunkard's walk on the unit sphere with step size
$\theta$ and show that the walk converges in order $C/\sin^2 \theta$
steps in the discrepancy metric ($C$ a constant).  
This is an application of techniques we develop for 
bounding the discrepancy of random walks on Gelfand pairs
generated by bi-invariant measures.
In such cases, Fourier analysis on the acting group
admits tractable computations involving spherical functions.
We advocate the use of discrepancy as a metric on probabilities 
for state spaces with isometric group actions.
\end{abstract}

\maketitle

\section{Introduction}
Fix $\theta \in (0,\pi)$. Consider the following random walk on
the unit sphere $S^2$ in $\R^3$, whose steps are geodesic arcs of
length $\theta$. (Such arcs subtend an angle of $\theta$ at the
center of the sphere). The random walk starts at the north pole,
and at each step a uniformly random direction is chosen and the
walk moves a geodesic distance $\theta$ in that direction.  We
refer to this walk as the {\em drunkard's walk on the sphere}.

The purpose of this paper is to develop techniques for bounding
the discrepancy metric for random walks on Gelfand pairs, using
the drunkard's walk as an example.  Our bounds are sharp enough to
give a rate of convergence.  Let $D(k)$ denote the discrepancy
distance (defined later) between the $k$-th step probability
distribution of the drunkard's walk and the uniform
(rotation-invariant) measure on $S^2$.  We show the following:

\begin{thm}
\label{maintheorem} For the drunkard's walk on the unit sphere
$S^2$ with step size $\theta$, the discrepancy of the walk after $k$
steps satisfies, for $k=\frac{C}{\sinsq}$,
$$
0.4330 \ e^{-C/2} \leq D(k) \leq 4.442 \ e^{-C/8}.
$$
\end{thm}

Thus order $C/\sinsq$ steps are both necessary and sufficient to
make the discrepancy distance of this walk from its limiting
distribution uniformly small.  The result makes intuitive sense,
since the number of steps to random should be large when $\theta$
is close to $0$ or $\pi$, and small when $\theta$ is close to
$\pi/2$.  Moreover, if $\theta \approx 1/n$ for large $n$, then
this result shows that order $n^2$ steps are necessary and
sufficient; this is similar to nearest-neighbor random walks on
$\Z/n\Z$ (e.g., see \cite{d}).  We also note that given $\theta$ this
walk does not exhibit a sharp cutoff phenomenon.

We frame our analysis in the context of a random walk on a
homogeneous space, i.e., a space with a transitive group action.
$S^2$ is a homogeneous space by the action of $SO(3)$. Although
the drunkard's walk is not generated by a group action, we show
its equivalence with a walk that is.

For random walks on groups, Fourier analysis is often used to
obtain rates of convergence.  On homogeneous spaces, we can lift
the walk to acting group and do Fourier analysis there, although
for non-commutative groups the group representations can be quite
complicated.  However, when the homogeneous space is a {\em
Gelfand pair} (as in this case), Fourier transforms of
bi-invariant measures and functions on the group simplify greatly,
allowing for tractable computations involving the spherical
functions.  In our example, the generating measure on $SO(3)$ can
be made bi-invariant, and the spherical functions are Legendre
polynomials.

Much of the literature on random walks on Gelfand pairs is
limited to discrete homogeneous spaces.  Diaconis \cite{d}
presents a survey and an annotated bibliography; applications
include walks on subspaces of vector spaces over finite fields
\cite{gr} and walks on $r$-sets of an $n$-set \cite{ds}. Rates
are given in the total variation metric.

Continuous examples have been addressed by Voit, who studied
families of isotropic random walks on spheres \cite{voit-sphere,
voit-ann} and other homogeneous spaces \cite{voit-twopt}.  Central
limit theorems are obtained using convergence in distribution or
total variation as the dimension $n\rightarrow \infty$.  Such
results differ from ours in that: (1) we work with a specific walk
rather than a family, e.g., we obtain explicit bounds for a
specific $n$ rather than asymptotic results for large $n$, (2) we
focus on rates of convergence of the walk on the homogeneous
space, rather than convergence of a central limit theorem on the
double coset space as in \cite{voit-sphere,voit-twopt}, and (3) we
use the discrepancy metric to measure convergence. We argue that
it is a natural metric to use for walks on homogeneous spaces, and
develop techniques to bound it. While we only illustrate our
methods on the $2$-sphere, similar methods can be used to give
explicit discrepancy bounds for walks on high-dimensional spheres
and other Gelfand pairs.

For walks on continuous groups, we mention the work of Rosenthal
\cite{rosen1} and Porod \cite{porod, porod-ann}, who obtain total
variation rates of convergence for random walks on $SO(n)$ and
other compact groups where the generating measures are
conjugate-invariant; this is another situation where the
representations simplify enough to get Fourier bounds.

This paper is organized as follows.  Section 2 gives background
on the discrepancy metric and justifies its use over other common
metrics on probabilities.  Section 3 develops several equivalent
formulations for the drunkard's walk.  Section 4 gives a
formulation with a bi-invariant generating measure.  This
simplifies the Fourier analysis in Section 5, where matching
upper and lower bounds for the convergence rate of the drunkard's
walk are derived (Theorems \ref{drunkardub} and \ref{drunkardlb}).
We summarize our methods for handling random walks on arbitrary
Gelfand pairs in Section 6.

For the uninitiated, Appendix A collects relevant background on
Fourier analysis on groups, Gelfand pairs, and representations of
$SO(3)$ that are needed to make this paper self-contained.
Appendix B contains proofs of technical results that are not
central to the development of the ideas in this paper.

\section{The discrepancy metric}

Let $(X,d)$ be a metric space with metric $d$.  Given any two
probability measures $P,Q$ on $X$, define the {\em discrepancy
distance} between $P$ and $Q$ by: $$ D(P,Q) = \sup_{\text{all
balls B}} |P(B) - Q(B)| $$ where a ``ball'' in $X$ denotes any
subset of the form $\{ x: d(x,x_0) \leq r\}$ for some $x_0 \in X$
and real number $r \geq 0$.  It is easy to check that the
discrepancy is a metric on probability measures.

When $X$ is the unit cube in $\R^n$ and $Q$ is Lebesgue measure on
$X$, this definition reduces to the notion of discrepancy commonly
used by number theorists to study {\em uniform distribution of
sequences} in the unit cube (e.g., see
\cite{drmotatichy, kn}).  Diaconis \cite{d} was perhaps the
first to suggest the use of discrepancy to measure rates of
convergence of random walks.  Su \cite{suthesis, su-trans,
su-leveque} explored properties of this metric and obtained sharp
rates of convergence for certain random walks on the hypercube,
circle, and torus.

We shall be concerned with the case where $X=S^2$ and the
metric on $S^2$ is inherited from its inclusion in $\R^3$. Thus
balls may be visualized as spherical ``caps'' on the sphere.  As
noted later, the group of rotations $SO(3)$ acts on $S^2$ in a
natural way and the metric on $S^2$ is invariant under this
action.  Thus images of balls under this action are still balls,
so the discrepancy metric on measures inherits this rotation
invariance.

Unlike the total variation metric (which is in more frequent use
among probabilists), the discrepancy metric recognizes both the
topology and the group action on the underlying space.  For
infinite compact state spaces this can be important.  For
instance, if $P$ is a probability measure on $S^2$ supported on a
finite set of points, and $Q$ is  the uniform (rotation-invariant)
probability measure, then the total variation distance between $P$
and $Q$ remains equal to 1 no matter how the points are arranged.
On the other hand, the discrepancy $D(P,Q)$ will capture how
``well-distributed'' the points in $P$ are.  Another example is a
simple random walk on the circle generated by an irrational
rotation (see \cite{su-trans}), which converges weak-* to Haar
measure; discrepancy captures this convergence, but total
variation is blind to it.  Thus the discrepancy metric is
well-suited to studying random walks on continuous state spaces
generated by isometric group actions.  

We favor the use of discrepancy over other common metrics (e.g.,
the Prokhorov, Wasserstein metrics) because there are tractable
bounding techniques for discrepancy involving Fourier
coefficients.  In fact, one of the main goals of this paper is to
show that we can develop upper and lower bounds for discrepancy
which give sharp rates of convergence in many cases because the
dominant terms in each expression match.

We remark that discrepancy bounds can be used to bound other metrics
by exploiting known relationships between them \cite{gibbs-su}; for instance,
the discrepancy is bounded above by total variation, so
discrepancy lower bounds also offer a way to obtain total variation
lower bounds.

A property that will be needed later is that on groups,
discrepancy decreases with convolution:
\begin{thm}
\label{convolvedecrease} If $P,Q,\nu$ are arbitrary probability
measures on a compact group $G$, then $$ D(P*\nu,Q*\nu) \ \leq \
D(P,Q) . $$ Hence when $Q=U$, the uniform (Haar) measure, we have
\begin{equation}
D(P*\nu,U) \ \leq \ D(P,U).
\end{equation}
\end{thm}
See \cite{su-trans} for a proof.

\section{The Drunkard's Walk and Equivalent Formulations}
The drunkard's walk on the sphere is not a random walk generated
by a group action.  However, we show in this section that it is
equivalent to one that is, in the sense that the two random walks
generate the same $k$-th step probability distribution even
though their observed behaviors may appear quite different.

Readers familiar with {\em hypergroups} may not be surprised by
the equivalence and the ensuing analysis, since the associated
double coset space of this walk is a commutative hypergroup, and
much of our analysis can be framed in that language.  We have
avoided it; interested parties are referred to \cite{bloom-heyer}.

Let $N$ be the isotropy subgroup of $SO(3)$ fixing $\np$, the
north pole. Let $E\subset SO(3)$ denote the set of all rotations
which fix a point on the equator and move the north pole by
geodesic distance $\theta$ along the surface of the sphere. Let
$Q$ denote the probability distribution supported on $E$ that is
left $N$-invariant.

\smallskip
{\bf Formulation 1.  The Drunkard's Walk.} This is the
walk considered at the opening of this paper; a drunkard starts at
the north pole and at each step picks a uniformly random direction
and advances along the sphere in that direction by geodesic
distance $\theta$.

\smallskip
 Let $\Y_k$ for $k=0,1,2,...$ denote random variables which
describe the location of the drunkard at time $k$.  Thus
$\Y_0=\np$.  If $g_i$ is an $SO(3)$-valued random variable with
values in $E$ and distribution $Q$, the position of the drunkard
at time $k$ is given by $\Y_k= g_1 g_2 \cdots g_k \np$.  Thus,
this walk is not a random walk in which the next position is
generated by applying group actions to the current position in
the walk.  However, the following random walk is:

\smallskip
{\bf Formulation 2.  The Potted Plant.} Consider a
potted plant initially at the north pole.  At each step, a
rotation is chosen randomly from $E$ according to $Q$ and
performed on the sphere. Thus the point currently over the north
pole moves a distance $\theta$ in any direction.  This induces a
motion of the potted plant, wherever it currently is.

\smallskip
Note that with the given generating set $E$, the potted plant is
moved a geodesic distance less than or equal to $\theta$ at every
step, since for each rotation in $E$, $\np$ is on the equator of
the rotation axis and hence moves the farthest.

If $g_i$ is a $SO(3)$-valued random variable with distribution
$Q$, the position of the potted plant at time $k$ is given by
$\Y_k = g_k g_{k-1} \cdots g_1 \np$. Since the $g_i$ are
independent and identically distributed, this shows that
Formulations 1 and 2 are equivalent and generate the same $k$-th
step probability distribution on the sphere.  This may be
surprising in light of the fact that the steps of the random walk
in Formulation 2 are smaller than in Formulation 1.

The next random walk, while not essential in what follows, also
generates the same $k$-th step probability distribution and we
mention it for the sake of interest.

\smallskip
{\bf Formulation 3.  Rotate and Spin.}  Fix any rotation
$R_\theta$ which displaces the current north pole by geodesic
distance $\theta$.  Start the random walk at the north pole, and
at each step perform $R_\theta$ followed by a uniform spin around
the north-south axis.  (The uniform spin moves the random walk to
a random point anywhere on the same latitude.)

\smallskip
Though $R_\theta$ is not necessarily contained in the
set $E$ defined earlier, it does yield the same $k$-th step
probability distribution as the previous formulations.  This may
be seen as follows.

Consider the double coset space $SO(3)//N$.  Each double coset is
characterized by the latitude to which it sends the north pole.
Thus $R_\theta=n' g_0 n''$ for some $n',n''\in N$ and a $g_0\in
E$.  Then $E=N g_0$.  Let $n_i$ denote an $N$-valued random
variable distributed according to Haar measure on $N$. The walk
description shows that at the $i$-th step, $n_i n' g_0 n''$ acts
on the random walk's current position. Therefore its position at
time $k$ is given by $ \Y_k = (n_k\ n' g_0 n'')( n_{k-1}\ n' g_0
n'')\cdots (n_1\ n' g_0 n'')\np$.  Since $n_i n' g_0$ and $n'' n_i
n' g_0$ are identically distributed according to $Q$ and $n''
\np=\np$, the above random variable has the same $k$-th step
distribution as the other formulations above.

\smallskip Our original goal was to study Formulation 1, the
drunkards' walk. Via the above equivalence we choose instead to
study Formulation 2, because it is a random walk generated by a
group action. However, the generating measure $Q$, while left
$N$-invariant, is not bi-invariant.
In light of Theorem \ref{gelfandequiv}, a 
bi-invariant generating measure would greatly simplify the ensuing
Fourier analysis.  (See Appendix A for background material on Fourier
analysis on compact groups and bi-invariant measures).
In the next section, we remedy
this problem by introducing a fourth random walk (Formulation 4) which is
equivalent to Formulation 2 and whose generating measure is
bi-invariant.

\section{A Bi-invariant Formulation}

We are interested in the discrepancy distance between the $k$-th
step distribution of the drunkard's walk and $U_{S^2}$, the
uniform (rotation-invariant) distribution on $S^2$. To simplify
notation, we write
\begin{equation}
\label{discreplaw} D(k) =D(\law (\Y_k),U_{S^2})
\end{equation}
where $\law (\Y_k)$ denotes the distribution of the random
variable $\Y_k$ in Formulation 2. We investigate the behavior of
$D(k)$ as a function of the number of steps $k$.

Recall that the homogeneous space $S^2$ can be regarded as the
left cosets of $N$ in $SO(3)$, so that the quotient map
$SO(3)\rightarrow S^2$ sends a rotation $g$ to the point $g\np \in
S^2$. A random walk on $S^2$ generated by an $SO(3)$-action (such
as Formulation 2) may then be regarded as a random walk
``upstairs'' on $SO(3)$ with an initial distribution $U_N$, Haar
measure on $N$ (which is the pre-image of the starting point
$\np$). The probability distribution upstairs evolves as usual
for a random walk on a group, so that after one step the
distribution is given by $Q*U_N$ and after $k$ steps by
$Q^{*k}*U_N$.  The probability of finding the original walk in a
ball $B \subset S^2$ is the same as finding the lifted walk on
$SO(3)$ in $\tilde{B}=BN \subset SO(3)$.  Hence
\begin{equation}
\label{discrepupstairs} D(k) = \sup_{\tilde{B}} |Q^{*k}*U_N
\,(\tilde{B})
     - U(\tilde{B})|.
\end{equation}
where $U$ is Haar measure on $SO(3)$ and the supremum is taken
over all $\tilde{B}$, pre-images of balls under the quotient map
$SO(3)\rightarrow S^2$.

At this point we would appeal to Fourier analysis to deal with the
convolutions above.  However, $Q$ is left $N$-invariant but not
bi-invariant; 
recall that 
we desire bi-invariance to simplify the Fourier analysis.

The following proposition shows that for random walks on groups,
averaging the generating measure $Q$ to make it bi-invariant will
affect the rate of convergence in discrepancy by at most one
step.  This result is the analogue of a result of Greenhalgh
\cite{gr}, who obtained a similar result for the total variation
distance.

\begin{prop}
\label{greenhalgh} Let $Q$ denote any left $N$-invariant
probability measure on a group $G$, let $U$ and $U_N$ denote Haar
measure on $G$ and $N$ respectively.  If $ \bar{Q}= Q*U_N$, then
$\bar{Q}$ is $N$-bi-invariant and $$ D(\bar{Q}^{*k},U) \leq
D(Q^{*k},U) \leq D(\bar{Q}^{*(k-1)},U). $$
\end{prop}

\begin{proof}
Left invariance for $Q$ means $U_N * Q = Q$.  We use this to
establish bi-invariance for $\bar{Q}$, which means $U_N * \bar{Q}
* U_N = \bar{Q}$. This follows from $U_N*\bar{Q}*U_N = U_N * (Q *
U_N) * U_N = (U_N
* Q)*(U_N*U_N) = Q *U_N = \bar{Q}$.

For the second assertion, note that
\begin{eqnarray*}
\bar{Q}^{*k} 
&=& (Q * U_N)^{*k} 
= Q * (U_N * Q)^{*(k-1)} * U_N 
= Q^{*k} * U_N .\\ 
Q^{*k} 
&=& Q * (U_N * Q)^{*(k-1)} 
= (Q * U_N)^{*(k-1)} * Q  
= \bar{Q}^{*(k-1)} * Q .
\end{eqnarray*}
Using Theorem \ref{convolvedecrease} we obtain $ D(Q^{*k}*U_N,U)
\leq D(Q^{*k},U)$ and $D(\bar{Q}^{*(k-1)}*Q,U) \leq
D(\bar{Q}^{*(k-1)},U)$, which with the above equations yield the
desired conclusion.
\end{proof}

Thus a random walk on a group with generating measure $Q$ differs
by no more than one step from a random walk proceeding according
to $\bar{Q}$, which may be viewed as the average of the measure
$Q$ over the left cosets of $N$.

However, for a random walk on a homogeneous space, even more can
be said if the walk begins at the point fixed by the isotropy
subgroup:

\begin{prop}
\label{drunkequiv4} Suppose $X$ be a homogeneous $G$-space with
isotropy subgroup $N$ fixing $x_0\in X$, and $Q$ is a
left-invariant probability on $G$.  Let $D(k)$ denote the
discrepancy of the random walk starting at $x_0$ and evolving via
a group action with elements chosen according to $Q$.  Let
$\bar{D}(k)$ denote the discrepancy of the random walk starting at
$x_0$, but evolving according to $\bar{Q}$. Then
$$\bar{D}(k)=D(k). $$
\end{prop}

\begin{proof}
This follows from the fact shown in the previous proof, that
$$
\bar{Q}^{*k}=Q^{*k} * U_N.
$$
The right side, when regarded as a measure on $S^2$, describes the
location of the $Q$-generated walk.  But by the right invariance
of $\bar{Q}$, the left side is equal to $\bar{Q}^{*k}*U_N$, which
when regarded as a measure on $S^2$, describes the location of
the $\bar{Q}$-generated walk.
\end{proof}

This shows that the following is equivalent to Formulation 2.

\medskip
{\bf Formulation 4.}  Let $\bar{E}$ denote the set of {\em all}
rotations in $SO(3)$ which move the north pole $\np$ by a fixed
angle $\theta$.  Let $\bar{Q}=Q * U_N$ be the bi-invariant
generating measure obtained by averaging $Q$ from Formulation 2.
Consider a potted plant which starts at $\np$ and is moved
according to the following rule: at each step, a rotation is
chosen randomly from $\bar{E}$ according to $\bar{Q}$ and
performed on the sphere. This induces a motion of the potted plant
to a new location.

\medskip
Observe that we are able to throw extra rotations in ``for free''
and still obtain the same $k$-th step probability distribution.
This may be surprising because with the extra generating elements
the step size of the potted plant is no longer bounded by
$\theta$, as it was in Formulation 2.  In fact, the potted plant
could be moved around rather wildly at each step.

\medskip
Exploiting this equivalence, we shall, in the sequel, work with
Formulation 4.  To save notation we write $Q$ for the bi-invariant
measure $\bar{Q}$.  Right-invariance for $Q$ yields $Q=Q*U_N$, which when
substituted into (\ref{discrepupstairs}) gives
\begin{equation}
\label{discrepupstairs2} D(k) = \sup_{\tilde{B}}
|Q^{*k}(\tilde{B})
     - U(\tilde{B})|.
\end{equation}
where the supremum is taken over all ball pre-images $\tilde{B}$.
Hence, the discrepancy $D(k)$ as defined in (\ref{discreplaw})
can now be analyzed using expression (\ref{discrepupstairs2}).

\section{A Rate of Convergence}
We now proceed to derive a rate of convergence for the drunkard's
walk on the sphere.  Several calculations require the facts
reviewed in Appendices A and B; we alert the
reader with references.

Let $B_{y,r}$ denote a ball of geodesic radius $r$ centered at
$y\in S^2$. Such balls look like spherical ``caps'' on $S^2$. Let
$\bt_{y,r}$ denote its pre-image ``upstairs'' in $SO(3)$.  To
reduce notation, write $\bt_r=\tilde{B}_{\np,r}$ for the pre-image
of a ball centered around $\np$.  Let $\delta_r$ denote the
indicator function of $\bt_r$ on $SO(3)$.

A key observation (see \cite{beck-chen}) for evaluating measures
on balls is they can be regarded as convolutions with indicator
functions on those balls, i.e., for any right-invariant measure
$\nu$ on $SO(3)$,
\begin{equation}
\label{beck-obs}
\nu(\bt_{y,r})= \nu *\delta_r(\yt)
\end{equation}
for any
$\tilde{y} \in yN$. This follows from $\nu *\delta_r(\yt) =
\int_{g\in\bt_r} d\nu(\yt g^{-1}) = \nu(\yt\cdot \bt_r)$.

From (\ref{discrepupstairs2}) we have
\begin{eqnarray}
\label{dqkusup} D(k)
 &=& \sup_{y,r} \left| Q^{*k}(\bt_{y,r})
     -U(\bt_{y,r}) \right|   \nonumber \\
&=& \sup_{y,r} \left| Q^{*k}*\delta_r\,(\yt)
     -U(\bt_{y,r}) \right| .
\end{eqnarray}

We wish to use Fourier inversion to derive bounds for these
expressions in terms of the Fourier coefficients.  We need
continuity of $Q^{*k} * \delta_r$ for $k\geq 2$:
\begin{prop}
\label{Qkoncap}
Let $Q$ be defined as in Formulation 4, and let
$\delta_r$ be denote the indicator function of $\bt_r$. Then
$Q^{*k} * \delta_r$ is continuous for $k\geq 2$.
\end{prop}
This is proved in Appendix B.

Hereafter, assume $k\geq 2$. We shall also assume for the moment
that $Q^{*k} * \delta_r$ has an absolutely convergent Fourier
series, which will be verified later in the course of our
computations.  Since $Q^{*k} * \delta_r$ is a continuous function
for $k\geq 2$, it is exactly equal to its Fourier series (Theorem
\ref{fourierinversion}), so that from 
(\ref{dqkusup}) and (\ref{fourierseries}) we have
\begin{eqnarray}
D(k)   &=& \sup_{y,r} \left| \sum_{n=1}^\infty (2n+1) \ Tr\left[
     \what{Q}^k(\rho_n) \
     \what{\delta_r}(\rho_n) \ \rho_n(\yt) \right]\right|
\nonumber \\
 \label{tracesum}
 &\leq& \sup_{y,r}
\sum_{n=1}^\infty (2n+1) \left| \ Tr\left[
     \what{Q}^k(\rho_n) \
     \what{\delta_r}(\rho_n) \ \rho_n(\yt) \right]\right|
\end{eqnarray}
where $\rho_n$ is the irreducible representation of $SO(3)$ of
dimension $(2n+1)$.  The trivial representation $\rho_0$ does not
appear here since it was cancelled in (\ref{dqkusup}) by
$U(\bt_{y,r})$.

\begin{remark}
\label{qdeltabiinv}
 Since $Q$ and $\delta_r$ are both
$N$-bi-invariant on $SO(3)$, by Theorem \ref{gelfandequiv} there
is a basis for the representations such that their transforms are
identically zero except in the $(1,1)$-th entry.  Any such basis
(e.g., the spherical harmonics) has its first basis element given
by the Legendre polynomials, which are the spherical functions
for the Gelfand pair $(SO(3),N)$.
\end{remark}

Hence $\rho_n(\yt)_{(1,1)}=P_n( \cos \gamma)$, where $\gamma$
is the geodesic distance of $y$ from $\np$. Since the product
of the transforms of $\what{Q}^k$ and $\delta_r$ are identically
zero except for the $(1,1)$-th element, the only diagonal element
changed by multiplication by $\rho_n(\yt)$ is the $(1,1)$-th
entry.  Hence the trace (\ref{tracesum}) reduces to
\begin{eqnarray}
D(k) &\leq&  \sup_{y,r} \sum_{n=1}^\infty (2n+1)
    \left|
     \what{Q}^k(\rho_n)_{(1,1)} \
     \what{\delta_r}(\rho_n)_{(1,1)}  \
     P_n(\cos \gamma) \right|  \nonumber \\
     \label{tracesum2}
&\leq&  \sup_r \sum_{n=1}^\infty (2n+1)
    \left|
     \what{Q}^k(\rho_n)_{(1,1)} \
     \what{\delta_r}(\rho_n)_{(1,1)}  \right|
\end{eqnarray}
where the second inequality follows from (\ref{lessthanone}).
Notice that the sum in (\ref{tracesum2}) is precisely the sum in
Theorem \ref{fourierinversion} that needs to be checked for
convergence in verifying that $Q^{*k} *\delta_r$ has an absolutely
convergent Fourier series.  Hence when we bound the above
expression we will also have validated our use of Fourier
inversion in our computations.

From (\ref{oneoneentry}), we have
\begin{equation}
\label{Qtransform} \what{Q}(\rho_n)_{(1,1)} =  P_n(\cos\theta)
\end{equation}
since $P_n$ is constant on the support of $Q$. Also, for a ball
$B_r$ of geodesic radius $r$ and $n\geq 1$, formula
(\ref{intondoublecoset}) gives
\begin{eqnarray}
\left| \what{\delta_r}(\rho_n)_{(1,1)} \right|
   &=&
   \left|  \half \int_{\cos r}^1 P_n(x) \ dx \ \right| \nonumber\\
   \label{deltatransform} &=&
        \frac{ | P_{n-1}(\cos r) - P_{n+1}(\cos r) | }{2 (2n+1) } \\
   &\leq &  \frac{1}{2n+1}.\label{deltatransform2}
\end{eqnarray}
The integral of $P_n$ follows from (\ref{legintegral}) and noting
that $P_n(1)=1$, and the inequality follows from
(\ref{lessthanone}).

Substitution of (\ref{Qtransform}) and (\ref{deltatransform2})
into (\ref{tracesum2}) yields
\begin{equation}
\label{legendresum} D(k) \leq \sum_{n=1}^\infty
|P_n^k(\cos\theta)| .
\end{equation}

To bound the Legendre polynomials, we use the following
well-known bound (see Jackson \cite[p.63]{j}):
\begin{prop}
\label{legendrebound1} For $P_n$, the $n$-th Legendre polynomial,
and any $\theta$, $$ |P_n(\cos\theta)|^2 \leq \frac{2}{\pi n
\sinsq}. $$
\end{prop}
We derive an alternate bound, suitable for small $\theta$:
\begin{prop}
\label{legendrebound2} For $P_n$, the $n$-th Legendre polynomial,
and $n\sinsq \leq .9$, $$ |P_n(\cos\theta)|^2 \leq 1-
\frac{n\sinsq}{4}. $$
\end{prop}

This bound is better than Proposition \ref{legendrebound1} when
$n\sin^2\theta < 2-\sqrt{4-\frac{8}{\pi}} \approx .794$.  It is
proved in Appendix B. Using Propositions \ref{legendrebound1} and
\ref{legendrebound2} and the bound $1-x\leq e^{-x}$, the sum in
(\ref{legendresum}) can be estimated:

\begin{eqnarray*}
\sum_{n=1}^\infty  | P_n^k(x) |
  &\leq& \sum_{n\leq B/\sinsq} e^{-n k\sinsq/8} \quad
+ \sum_{n > B/\sinsq}
     \left(\frac{2}{\pi n \sinsq}\right)^{k/2} \\
&\leq& \frac{e^{-k\sinsq /8}}{1 - e^{-k\sinsq /8}}
      \ \ +\ \ \left(\frac{2}{\pi \sinsq}\right)^{k/2}
        \sum_{n > B/\sinsq}\frac{1}{n^{k / 2}}
\end{eqnarray*}
where $B=.9$.  Note that
\begin{eqnarray*}
\sum_{n > B/\sinsq}\frac{1}{n^{k / 2}} &\leq& \quad
\int_\frac{B}{\sinsq}^\infty \frac{dx}{x^{k/2}}
  \quad + \quad \left( \frac{\sinsq}{B} \right)^\frac{k}{2}\\
&=& \frac{2}{k-2}\left(\frac{\sinsq}{ B}\right)^{\frac{k}{2}-1}
          + \left( \frac{\sinsq}{B} \right)^\frac{k}{2}.
\end{eqnarray*}
Thus
\begin{eqnarray*}
 \sum_{n=1}^\infty | P_n^k(x) |
   &\leq& \frac{e^{-k\sinsq /8} }{1 - e^{-k\sinsq /8}}
      \ \ + \ \ \left(\frac{2}{\pi B}\right)^{k/2}
      \left( \frac{2 B }{(k-2) \sinsq}   +1 \right).
\end{eqnarray*}
Note that $\left(\frac{2}{\pi B}\right)^{1/2} <
e^{-1/8} < e^{-\sinsq/8}$. For $k=\frac{C}{\sinsq}$ and $C\geq 4$,
one sees that $(k-2)\sinsq \geq 2$ and $k\geq 4$, so that
\begin{eqnarray*}
\sum_{n=1}^\infty  | P_n^k(x) | &\leq& e^{-k\sinsq/8} \left(
\frac{1}{1-e^{-1/2}} +
    B + 1  \right)\\
&\leq& 4.442 \ e^{-C/8}.
\end{eqnarray*}
The above bound, together with (\ref{legendresum}), 
proves the following theorem.  (Note that the $C \geq 4$ restriction above is 
not needed below because the discrepancy $D(k)$ never exceeds 1.)
\begin{thm}
\label{drunkardub} For the drunkard's walk on the sphere with step
size $\theta$, the discrepancy after $k$ steps satisfies, for
$k=\frac{C}{\sinsq}$, $$ D(k) \leq 4.442 \ e^{-C/8}.
$$
\end{thm}
Thus order $\frac{C}{\sinsq}$ steps are sufficient to make the
discrepancy uniformly small.  The following lower bound confirms
the order is correct.
\begin{thm}
\label{drunkardlb} For the drunkard's walk on the sphere with step
size $\theta$, the discrepancy after $k$ steps satisfies, for
$k\geq 2$,
$$ D(k) \geq \frac{\sqrt{3}}{4} |\cos \theta|^k .
$$ For $k=\frac{C}{\sinsq}$, we have $$ D(k) \geq 0.4330\
e^{-C/2}. $$
\end{thm}
Thus order $\frac{C}{\sinsq}$ steps are needed to make the
discrepancy distance uniformly small.  Together, Theorems
\ref{drunkardub} and \ref{drunkardlb} prove Theorem
\ref{maintheorem}.

One way to obtain a lower bound for discrepancy is to evaluate the
difference of $Q^{*k}$ and $U$ on well-chosen ball.  The same
idea can be used for the total variation; one way to choose such
a ball (see \cite[p.29]{d}) is to take a set cut out by a random
variable consisting of the dominant terms in the Fourier series
of $Q^{*k}$. The mean and variance of the random variable and an
appeal to Chebyshev's inequality yield an estimate for $Q^{*k}$
on that set.

However, the proof of Theorem \ref{drunkardlb} illustrates a
different approach using ideas similar to those used in
\cite{su-leveque} for bounds on the torus.  We construct a ``local
discrepancy'' function which at each point evaluates the
discrepancy of the measure on a set of geodesic radius $r$
centered at that point. The function is bounded above by the total
discrepancy. As before, it can be rewritten in terms of a
convolution of the original measure and the indicator function of
the set.  An appeal to Plancherel's identity gives a sum with only
non-negative terms, so the dominant term can be pulled out as a
lower bound for discrepancy.

We remark that since discrepancy is a lower bound for total
variation, this lower bounding technique can also be used to
obtain lower bounds for random walks under total variation.

\begin{proof}
Define, for $g\in SO(3)$, $$ \Delta_r(g) = Q^{*k}(\tilde B_{y,r})
- U(\tilde B_{y,r}) $$ where $y$ is the image of $g$ under the
quotient map from $SO(3)$ to $S^2$. From (\ref{discrepupstairs2}),
we see that $\Delta_r(x) \leq D(k)$, and hence for all $r$,
\begin{equation}
\label{anyr} \int_{SO(3)} \Delta_r^2(g) \ d\mu \leq D(k)^2.
\end{equation}
On the other hand, Plancherel's identity on $SO(3)$ \cite[p.256]{dym-mckean}
yields
\begin{equation}
\label{plancherel} \int_{SO(3)} \Delta_r^2(g) \ d\mu \ = \
\sum_{n=0}^{\infty} (2n+1) \ Tr\left[ \what\Delta_r(\rho_n)
\what\Delta_r(\rho_n)^* \right]
\end{equation}
where the $^*$ denotes the conjugate transpose (here only). Notice
that $\Delta_r$ may be rewritten as:
$$ \Delta_r(x) = Q^{*k} * \delta_{\tilde B_r} (x) - U *
\delta_{\tilde B_r} (x) = (Q^{*k} - U) * \delta_{\tilde B_r} (x)
.
$$
Then $\what\Delta_r$ may be computed as $\what\Delta_r
(\rho_n)= (\ \what Q^k (\rho_n) - \what U(\rho_n) \ ) \
\what\delta_{\tilde B_r} (\rho_n)$.

For $n= 0$,
$ \what\Delta_r (\rho_0) = 0 $ since
$ \what Q^k(\rho_0) - \what U (\rho_0) = 1 - 1 = 0$.

For $n\neq 0 $, a trivial computation shows $\what U(\rho_n) = 0$,
and thus
$$
\what\Delta_r (\rho_n) = \what Q^k (\rho_n)\  \what \delta_{\tilde B_r}
(\rho_n).
$$
Remark \ref{qdeltabiinv} and the computations from Equations (\ref{Qtransform})
and (\ref{deltatransform}), when substituted into (\ref{plancherel}), and
combined with (\ref{anyr}), give
\begin{equation}
\label{drunkardsum} D(k)^2 \ \geq \ \sum_{n=1}^{\infty} (2n+1)
\left| P_n^k(\cos\theta) \left( \frac{P_{n-1}(\cos r) -
P_{n+1}(\cos r)}{2(2n+1)}\right) \right|^2
\end{equation}
where $r$ may be chosen arbitrarily. Taking only the dominant term
($n=1$) in the above expression, and letting $r=\pi/2$, we have
$\cos r=0$, $P_0(0)=1$, $P_2(0)=-1/2$, and $P_1(x)=x$. It follows
that $$ D(k) \geq \frac{\sqrt{3}}{4} |\cos\theta|^k $$ as was to
be shown. The second inequality in the theorem follows from $|\cos
\theta|^k = e^{\frac{k}{2} \ln \cos^2\theta }\geq
e^{-\frac{k}{2}\sin^2 \theta} $, using the fact that $\ln (1-x)
\geq -x$ for all $x$.
\end{proof}

For a tighter lower bound, one may use more terms in
(\ref{drunkardsum}) and adjust the choice of $r$; however, the
dominant term sufficed to obtain matching upper and lower bounds
for this random walk.

\section{Conclusion}
A similar analysis can be carried out for the discrepancy
convergence of any random walk on a Gelfand pair, when the
spherical functions are known.  Proposition \ref{drunkequiv4}
shows that making a generating measure bi-invariant will not
affect the rate of convergence.  The upper bound is obtained via
(\ref{beck-obs}), Fourier inversion to yield (\ref{legendresum}), 
and bounds on the
appropriate spherical function (e.g., Prop.
\ref{legendrebound2}).  The lower bound is obtained via Plancherel's
identity applied to the square of the local discrepancy function,
e.g., equations (\ref{anyr}) and (\ref{plancherel}), 
then choosing as many terms as needed.

We remark that our pair of strategies often works well for
obtaining matching upper and lower bounds because if there is a
dominant Fourier coefficient, it appears to the same order in
both upper and lower bounds.  In our example, $\what{Q}(1)$ was
the dominant term; compare the upper bound (\ref{legendresum})
and lower bound (\ref{drunkardsum}).  See 
\cite{hensley-su, su-trans, su-leveque} 
for more examples of this phenomenon in discrepancy bounds for 
random walks on groups.

\section*{Appendix A.}
This appendix reviews material on harmonic analysis, homogeneous
spaces, Gelfand pairs, representations of $SO(3)$, and Legendre
polynomials.

\subsection*{Fourier Analysis on a Compact Group}
A standard reference is the encyclopedic account by Hewitt and Ross
\cite{hr}.  Diaconis \cite{d} gives a concise introduction to Fourier
analysis on finite groups.  Dym and McKean \cite{dym-mckean} is a
readable introduction to Fourier series on $SO(3)$.

We assume henceforth that all compact groups are separable and
metrizable.
For any compact group $G$ there is a unique measure $\mu$ on $G$,
called (normalized) {\em Haar measure}, such that $\mu$ is
$G$-invariant and $\mu(G)=1$.

Let $V$ be a finite dimensional vector space over $\C$, 
the complex numbers.
Recall that a {\em representation} of a group $G$ on $V$
is a homomorphism $\rho:G\rightarrow GL(V)$.  
If $V$ has dimension $n$, then $\rho$ is said to have
{\em dimension} $n$.
A basis for $V$ can be chosen so that the image of $\rho$ 
with respect to this basis are unitary matrices.
If there is no non-trivial subspace of $V$ invariant under the action of
$G$, then $\rho$ is said to be {\em irreducible}; otherwise 
$\rho$ decomposes as a direct sum of irreducible representations.
(One can similarly define a representation $\rho$ of $G$ on a Hilbert
space, though if $G$ is compact $\rho$ decomposes into a direct sum of 
unitary representations of finite dimension.)

Two representations
$\rho$ on $V$ and $\rho'$ on $V'$ are {\em equivalent} 
if there is an isomorphism $\tau:V \rightarrow V'$ such that 
$\tau \circ \rho = \rho' \circ \tau$.  Let $\Sigma$ denote a
the set of equivalence classes of irreducible representations of $G$.
For a compact group, 
$\Sigma$ is countable and
furthermore, all the irreducible representations are finite
dimensional.

\begin{defn}
The {\em Fourier transform} of a complex-valued function $f$ on a
compact group $G$ at a representation $\rho$ of $G$ 
is defined by
$$
\what{f}(\rho) = \int_{g\in G}  f(g)\ \rho(g^{-1}) \ d\mu
$$
where $\mu$ is Haar measure on $G$.

Similarly, the Fourier transform of a measure $\nu$ on $G$ at
$\rho$ is defined by
$$
\what{\nu}(\rho) = \int_{g\in G} \rho(g^{-1}) \ d\nu(g).
$$
\end{defn}

We show how a function may be recovered from its Fourier transforms at
irreducible representations.  
Let $d_{\rho}$ denote the dimension of a representation 
$\rho$.
For any operator $A$, let $Tr[ A ]$ denote the trace of $A$, and 
let $\|A\|_{\varphi_1}$
denote the sum of the eigenvalues of the operator square root of
$AA^*$.
(Here, $^*$ denotes conjugate transpose.) 

\begin{defn}
For any $f\in L^1(G,\mu)$, the series 
\begin{equation}
\label{fourierseries} \sum_{\rho \in \Sigma}
      d_{\rho} \ Tr\left[ \what{f}(\rho) \ \rho(g)\right]
\end{equation}
is called the {\em Fourier series} of $f$.
(There is mild abuse of notation here: by $\rho \in \Sigma$, 
we really mean to choose a representative $\rho$ from each
class of irreducible representations in $\Sigma$.)
If
\begin{equation}
\label{absconvcondition} \sum_{\rho \in \Sigma}
  d_{\rho} \ \| \what{f}(\rho) \|_{\varphi_1}< \infty \ ,
\end{equation}
then 
$f$ is said to have an {\em absolutely convergent Fourier series}
\cite[(34.4)]{hr}. 
\end{defn}

\begin{thm}[Fourier inversion]
\label{fourierinversion} 
If a function $f$ on $G$
has an absolutely convergent Fourier series, then the 
Fourier series of $f(g)$ converges uniformly to a continuous function 
$\bar f(g)$, and $f(g)= \bar f(g)$ almost everywhere on $G$
with respect to Haar measure $\mu$.
\end{thm}

\begin{proof}
This theorem is embedded in Hewitt and Ross \cite{hr}, but obscured by
their exotic notation.  We briefly indicate how to ``prove'' this
theorem from results cited in \cite{hr}.

The set of functions with absolutely convergent Fourier series is
denoted in \cite{hr} by a symbol 
that resembles $\mathcal{R}(G)$, defined in (34.4).
Theorem (34.6) in \cite{hr} shows that any $f\in \mathcal{R}(G)$ is
equal almost everywhere to its Fourier series.  
Theorem (34.5.ii) shows that this Fourier series converges uniformly to
a continuous function that we have denoted $\bar f$.
\end{proof}

We remark that since the notation in Hewitt and Ross \cite{hr} is 
cumbersome and tough to wade through, for the sake of probabilists
we have simplified it by following the notation of Diaconis \cite{d}.  
To aid the reader wishing to follow the results quoted above,
we provide a ``dictionary'' between the two sets of notation:
in Hewitt and Ross \cite[(27.3)]{hr}, 
$\sigma$ denotes a class of equivalent
irreducible representations in $\Sigma$ 
and $U$ is a representative of that class;
we avoid reference to $\sigma$ (to eliminate an unnecessary layer of notation)
and use $\rho$ instead of $U$.
Hewitt and Ross denote an arbitrary element
of a group $G$ by $x\in G$; we use $g \in G$.
Their notations $A_\sigma$ and $U^{(\sigma)}_x$
refer to operators that correspond to our $\what{f}(\rho)$ and
$\rho(g)$, respectively (see \cite[(34.1.i), (34.4.a)]{hr}).

Note that if $f$ is continuous, Theorem \ref{fourierinversion} 
implies that 
if $f$ has an absolutely convergent Fourier series, 
it equals its Fourier series at {\em every} point.

\subsection*{Homogeneous Spaces and Gelfand Pairs}
Diaconis \cite[Chap.\ 3F]{d} provides an introduction to Gelfand
pairs on finite groups and an annotated bibliography. Dieudonne
\cite{dieudonneshort} is a concise introduction to Gelfand pairs
on compact and locally compact groups.

\begin{defn}
Let $G$ be a compact group and $X$ be a topological space.  An
{\em action} of $G$ on $X$ is a continuous mapping from $G \times
X \rightarrow X$ denoted by $(s,x) \mapsto s\cdot x = sx$ such
that $id\cdot x=x$ and $s\cdot(t \cdot x)=(st)\cdot x$.

If $G$ acts transitively on $X$, that is, if for any $x,y\in X$
there exists an $s$ such that $sx=y$, we call $X$ a {\em
homogeneous space}.
\end{defn}

Given a point $x_0 \in X$, let $N$ denote the {\em isotropy}
subgroup of $G$ with respect to $x_0$, i.e., the set of group
elements which fix $x_0$.  By construction, $N$ is a closed
subset of $G$.  The canonical isomorphism of $X$ onto $G/N$, the
left cosets of $N$, respects the action of $G$. Thus $g: xN
\mapsto (gx)N$.

Let $\mu_X$ denote the $G$-invariant measure on $X$ induced by Haar
measure on $G$.  Let $L^2(X)$ denote the space of all
complex-valued square-integrable functions on $X$ with respect to
$\mu_X$. The action of $G$ on $X$ induces an action of $G$ on
$L^2(X)$ by $g \cdot f(x)= f(g^{-1}x)$. This action is a 1-to-1
linear mapping of the vector space $L^2(X)$ into itself and so
defines a representation of $G$.

\begin{defn}
\label{bi-inv}
A function $f$ on $G$ is said to be 
{\em $N$-bi-invariant} 
if $f(n' g n'')=f(g)$ for all $n',n'' \in N$ and
$g\in G$.  A measure $\nu$ on $G$ is {\em $N$-bi-invariant} if
for any measurable set $A$ in $G$, $\nu(n' A n'')=\nu(A)$ 
for all $n',n'' \in N$.  Thus $\nu$
satisfies $\nu * \mu = \mu * \nu = \nu$ where $\mu$ is Haar measure on $G$.
\end{defn}

In this paper bi-invariance on a homogeneous space $G$ 
will be understood to mean with
respect to the isotropy subgroup $N$.  Note that bi-invariant
functions on $G$ are constant on double cosets $NgN$ and may therefore
be viewed as functions on the double coset space (denoted $G//N$), 
or as left-invariant functions
on $X$ via its isomorphism with $G/N$.

\begin{defn}
The pair $(G,N)$ is called a {\em Gelfand pair} if the convolution
algebra $L^2(G//N)$ of $N$-bi-invariant functions is commutative.
We sometimes say $X\cong G/N$ is a Gelfand pair 
when $G$ is understood by context.
\end{defn}

The next fact about Gelfand pairs is the most important 
for our purposes.
A similar result for the finite group context may be found 
in \cite[p.54]{d}.

\begin{thm}
\label{gelfandequiv}
If $(G,N)$ is a Gelfand pair, then
for every irreducible representation $\rho:G\rightarrow
GL(V)$ there is a basis of $V$ such that for all 
functions $f$ (resp. measures $\nu$) bi-invariant with respect to $N$,
the Fourier transform
$\what{f}(\rho)$ (resp. $\what{\nu}(\rho)$) in that basis
contains only zeroes except
possibly for the $(1,1)$-th entry.
\end{thm}

\begin{proof}
Dieudonne \cite[(22.5.6)]{dieu-t} shows the algebra
$L^2(G//N)$ is commutative if and only if the number of times 
the trivial representation appears in $\rho|_N$, the
the restriction of $\rho$ to $N$, is zero or one.

If one, this trivial representation corresponds to
a one-dimensional subspace of $V$ fixed by $N$, i.e., the
left $N$-invariant functions on $X$; choose the unique function
$s(x)$ on $X$ normalized so that $s(x_0)=1$.  This is sometimes
called the {\em spherical function} of $(G,N)$ corresponding to the
representation $\rho$.  Complete $s$ to a basis for $V$ so that the 
matrices of $\rho|_N$ break into irreducible ``blocks''.  
Then for a right $N$-invariant function $f$:
\begin{eqnarray}
\what{f}(\rho) 
&=& \int_{g\in G}  f(g)\ \rho(g^{-1}) \ d\mu \nonumber \\
&=& \int_{n\in N} \int_{x\in G/N} 
	f(xn)\ \rho(n^{-1}x^{-1})\ d\mu_X\ d\mu_N \nonumber \\
&=& \int_{n\in N} \rho(n^{-1})\ d\mu_N \cdot
	\int_{x\in G/N} f(x) \rho(x^{-1}) \ d\mu_X \label{twointegs} \quad,
\end{eqnarray}
where $\mu_N$ is Haar measure on the subgroup $N$.  The second
equality is obtained by choosing a coset representative $x$ 
from each coset in $G/N$ and expressing $g=xn$ for some $x$ and $n\in
N$, and noting that Haar measure $\mu$ decomposes as a product measure
$\mu_X \cdot \mu_N$.
A similar argument holds for a right-invariant measure $\nu$, noting
that $\nu$ decomposes as a product measure $\nu_X \cdot \mu_N$ because
of right-invariance.

By the orthogonality relations for matrix entries
\cite[(21.2.5.c)]{dieu-t}
of irreducible representations of $N$, 
the left-most integral of (\ref{twointegs}) 
produces a matrix consisting of zeroes except possibly for the
$(1,1)$-th entry. 
Thus $\what{f}(\rho)$ (resp. $\what{\nu}(\rho)$) 
has zero entries except possibly for the first
row.  A similar argument using the left-invariance of $f$ (resp. $\nu$) 
shows that $\what{f}(\rho)$ (resp. $\what{\nu}(\rho)$) 
has zero entries except possibly for the first
column.  Together, these statements imply
that the only entry that could possibly be non-zero is the $(1,1)$-th entry.

If the trivial representation does not appear in the
the restriction of $\rho$ to $N$, the argument above holds by ignoring
the role of $s$ when choosing a basis for $V$.  Orthogonality then shows 
that the left-most integral of (\ref{twointegs}) yields a zero matrix.
\end{proof}

There is thus one spherical function $s_i(x)$ for every irreducible
representation $\rho_i$ appearing in $L^2(X)$.  These induce 
$N$-bi-invariant functions $\tilde{s}_i$ on $G$.
In the theorem above the $\tilde{s}_i(g)$ appears as the
$(1,1)$-th entry of $\rho_i(g)$ for an appropriate basis.  Hence
for any measurable function $f$ on $G$, the $(1,1)$-th entry of
the Fourier transform at $\rho_i$ satisfies
\begin{equation}
\label{foneoneentry} [\what{f}(\rho_i)]_{(1,1)}= \int_{g\in G}
f(g) \ \tilde{s}_i(g) \ d\mu.
\end{equation}
Similarly, for a measure $\nu$ on $G$,
\begin{equation}
\label{oneoneentry} \ [\what{\nu}(\rho_i)]_{(1,1)}= \int_{g\in G}
\ \tilde{s}_i(g) \ d\nu .
\end{equation}
Dieudonne \cite{dieudonneshort} is a readable introduction to the
general theory of spherical functions; Letac \cite{letac}
computes them in several examples.

\subsection*{The sphere as a Gelfand pair}
The rotation group $SO(3)$ acts on the unit sphere $S^2$ by the
natural inclusion of $S^2$ in $\R^3$.  This action is clearly
transitive on $S^2$, so $S^2$ is a homogeneous space. In fact
arises from the Gelfand pair $(G,N)$, where $G$ is the rotation
group $SO(3)$, and $N$ is the isotropy subgroup of rotations
fixing $\np$, the north pole. By restriction to the plane
orthogonal to $\np \in \R^3$, we see that $N$ is isomorphic to
the group $SO(2)$. The sphere $S^2$ may then be regarded as the
space $SO(3)/SO(2)$. In fact, for all $n\geq 2$, Dieudonne
\cite{dieudonneshort} shows that $S^n \cong SO(n+1)/SO(n)$ is a
Gelfand pair.

\subsection*{Representations of $SO(3)$}
For a good reference on representations of $SO(3)$ and other
compact Lie groups, see Brocker and tom Dieck
\cite{brockertomdieck}.  

Let $\Delta=\frac{\del^2}{\del x_1^2} + \frac{\del^2}{\del
x_2^2}+\frac{\del^2}{\del x_3^2}$ be the Laplace operator on
$\R^3$.  The {\em harmonic polynomials} are the set of all
complex-valued homogeneous polynomials $f$ in $x_1,x_2,x_3$ of
degree $n$ such that $\Delta f =0$; the restrictions of these
functions to the sphere $S^2$ form a set $V_n$, the {\em spherical
harmonics} of degree $n$ (one of which is the spherical function
$s_n$).

The action of $SO(3)$ on $V_n$ is induced by its action on $\R^3$
in the manner described earlier:  $g\cdot f(x) = f(g^{-1}x)$.
Moreover, $V_n$ is irreducible and finite-dimensional, and {\em
every} irreducible representation of $SO(3)$ arises in this way.
The dimension of $V_n$ is $2n+1$.

\subsection*{Legendre polynomials}
The spherical functions $s_i$ on $S^2$ are given by the well-known
Legendre polynomials $P_i$ in the following way:  for $y\in S^2$,
$s_i(y)=P_i(x)$ where $x=\cos \theta_y \in [-1,1]$ and $\theta_y$
is the geodesic distance between $y$ and $\np$ on $S^2$.  Just as
$S^2$ is (isomorphic to) the left cosets of $N$ in $G$, the set
$[-1,1]$ is the double coset space of this Gelfand pair.

Since Haar measure on $SO(3)$ induces the uniform
(rotation-invariant) probability measure on $S^2$ and uniform
probability measure on $[-1,1]$, we can compute
(\ref{foneoneentry}) as
\begin{equation}
\label{intondoublecoset} [\what{f}(\rho_n)]_{(1,1)} =
\int_{SO(3)} f(g) \tilde{s}_n(g) \ d\mu = \frac{1}{2} \int_{-1}^1
f(x) P_n(x) \ dx
\end{equation}
where $\mu$ denotes normalized Haar measure on $SO(3)$, $dy$
denotes the uniform measure on $S^2$, and $dx$ is Lebesgue
measure on $\R$.  See \cite[p.239]{dym-mckean}.

The Legendre polynomials have the generating function
\begin{equation*}
\sum_{i=0}^\infty P_i(x)\ r^i = \frac{1}{ (1-2xr +r^2)^{1/2}}
\end{equation*}
and the first few are: $P_0(x)=1, P_1(x)=x, P_2(x)=\half(3x^2-1),
P_3=\half(5x^3-3x)$. Furthermore, for all $n\geq 0$, $P_n(1)=1$
and
\begin{equation}
\label{lessthanone}
 |P_n(x)|\leq 1
\end{equation}
for $x \in [-1,1]$. The following identity will be needed later.
For $n\geq 1$,
\begin{equation}
\label{legintegral}
   P_n(x) =  \frac{ P'_{n+1}(x) - P'_{n-1}(x)}{2n+1}
\end{equation}
This follows from the generating function for $P_n(x)$.  A nice
account of these and other properties of Legendre polynomials may
be found in Jackson \cite{j}.

\section*{Appendix B.}
This appendix contains the proofs of some technical results
(Proposition \ref{Qkoncap} and Proposition \ref{legendrebound2})
that are not central to the development above.

\subsection*{Proof of Proposition \ref{Qkoncap}}
To show that $Q^{*k}
* \delta_r\,(x)$ is continuous for $k\geq 2$, we first require the
following technical lemma.

\begin{lemma}
\label{convolvecontin} Let $\nu$ be a positive measure on a
compact metric group $G$, and let $f$ be any measurable, bounded
function $f$ with discontinuities on a set $D_f$.  Let $xD_f^{-1}$
denote the set $\{xd^{-1}:d\in D_f\}$.  Given $x$, if $\nu(x
D_f^{-1})=0$, then the convolution $$ h(x)= \nu*f \, (x) = \int_G
f(z^{-1}x)\ d\nu(z) $$ is continuous at $x$.
\end{lemma}

\begin{proof}
To show $h(x)$ is continuous at $x$, consider any sequence $x_n
\in G$ such that $x_n\rightarrow x$.  It must be shown that
$h(x_n)\rightarrow h(x)$.

Let $w_n(z)=f(z^{-1}x_n)$ and $w(z)=f(z^{-1} x)$.  Since $f$ is
bounded, all the $w_n$ and $w$, being translates of $f$, are
uniformly bounded by some constant function.  This constant
function is in $L^1(G,\nu)$, since $G$ is compact.

Also, $w_n(z)\rightarrow w(z)$ pointwise for all $z\not\in x
D_f^{-1}$, since $f$ is continuous at those points.  By the
assumption on $D_f$ we have pointwise convergence almost
everywhere.  By Lebesgue's dominated convergence theorem, $\int_G
w_n(z)\ d\nu(z) \rightarrow \int_G w(z)\ d\nu(z)$, which is
precisely the statement $h(x_n)\rightarrow h(x)$.  This completes
the proof of Lemma \ref{convolvecontin}.
\end{proof}

We can now prove Proposition \ref{Qkoncap}.

\begin{proof}
We apply  Lemma \ref{convolvecontin} setting $f(x)=\delta_r(x)$
and $\nu=Q$.  By inspection, $\delta_r$ is bounded by $1$. The
lemma implies $Q* \delta_r(x)$ is continuous everywhere except
possibly at $x=id$.  This may be seen by observing that the
discontinuity set of $\delta_r$ is $\partial \bt_r$, the boundary
of $\bt_r$. This is the pre-image of a circle on $S^2$.  On the
other hand, $Q$ regarded as a measure on $S^2$ is uniformly
supported on a circle at latitude $\theta$ from the north pole.
Any two circles on $S^2$ intersect in at most two points, unless
they are identical.  Hence $Q(\partial \bt_r)=0$ unless the
support of $Q$ intersects $\partial \bt_r$, which only occurs when
$\bt_r=\bt_{\np,\theta}$. This corresponds to a discontinuity in
$Q*\delta_r(x)$ at $x=id$ when $r=\theta$.

We now apply Lemma \ref{convolvecontin} again to show that $
Q^{*k}*\delta_r(x)$ is continuous for $k= 2$.  The preceding
observations show that $Q* \delta_r(x)$ is continuous almost
everywhere (except possibly at the identity which is not in the
support of $Q$).  It is bounded by $1$.  Applying the lemma for
$f(x)=Q* \delta_r (x)$ and $\nu=Q$ shows that $Q^{*2}*\delta_r
(x)$ is continuous everywhere.

Now proceed by induction on $k$.  For $k\geq 3$, let $\nu=Q$ and
let $f(x)=Q^{*(k-1)}*\delta_r(x)$, which is continuous.  Then
Lemma \ref{convolvecontin} shows that $\nu * f(x)=Q^{*k}*\delta_r(x)$
is continuous.
\end{proof}

\subsection*{Proof of Proposition \ref{legendrebound2}}  This
proves a Legendre bound for small $\theta$.

\begin{proof}  Let $x=\cos\theta$.  From \cite[p.62]{j},
$$ |P_n(x)| \leq \frac{2}{\pi}\int_0^{\pi/2}e^{-n z^2 w^2 /2}\
dw $$ for $z=\frac{2}{\pi}(1-x^2)^{1/2}$.
  Substitute
$t=n^{1/2}z w$, and set $A=n^{1/2} z \pi / 2 = n^{1/2}\sin\theta$.
Obtain $$ |P_n(x)| \leq \frac{1}{A}\int_0^A e^{-t^2/2} \ dt. $$

To estimate this integral, square both sides and consider the
double integral over a square in first quadrant of the plane, and
then change to polar coordinates:
\begin{eqnarray*}
|P_n(x)|^2 &\leq& \frac{1}{A^2} \int_0^A\ \  \int_0^A
     e^{(-t_1^2 -t_2^2)/2}\ dt_1 \ \ dt_2 \\
&\leq& \frac{2}{A^2} \int_{\phi=0}^{\pi/4} \int_{r=0}^{A/\cos\phi}
     e^{-r^2/2} \ r \ dr \ d\phi \\
&\leq& \frac{2}{A^2} \int_{\phi=0}^{\pi/4} (1 -
e^{-A^2/2\cos^2\phi})\
  d\phi\ .
\end{eqnarray*}
For $0\leq\phi\leq\pi/4$, $y=\frac{A^2}{2}\sec^2\phi \leq A^2$
which by assumption is less than $.9$.  Now for $y< .9$, the
inequality $1-e^{-y} \leq y- \frac{3y^2}{8}$ holds, and yields
$$
|P_n(x)|^2 \leq \frac{2}{A^2} \int_{\phi=0}^{\pi/4}
    \left( \frac{A^2}{2}\sec^2 \phi - \frac{3A^4}{32}\sec^4\phi\right)
     \ d\phi \ .
$$
Integrating the right hand side gives
\begin{equation*}
|P_n(x)|^2 \leq
     \left. \left[
	\tan\phi 
     -  \frac{3A^2}{16} 
	\left(\frac{\tan^3\phi}{3}
     +\tan\phi\right) \right] \right|_0^{\pi/4}
     \leq 1 - \frac{A^2}{4}.
\end{equation*}
\end{proof}


\begin{thebibliography}{[DGM]}

\bibitem{beck-chen} Beck, J. and Chen, W. {\em Irregularities of
Distribution}, Cambridge Univ. Press, 1987.

\bibitem{bloom-heyer}  Bloom, W.R. and Heyer, H.  {\em Harmonic
Analysis of Probability Measures on Hypergroups}, deGruyter,
Berlin, 1994.

\bibitem{brockertomdieck}  Brocker, T. and tom Dieck, T.
{\em Representations of Compact Lie Groups}.  Springer-Verlag, 1985.

\bibitem{d}  Diaconis, P.  {\em Group Representations in Probability
and Statistics.}  IMS Lecture Notes-Monograph Series, vol. 11, 1988.

\bibitem{ds}  Diaconis, P. and Shahshahani, M.
Time to reach stationarity in the Bernoulli-Laplace diffusion
model. {\em SIAM J. Math. Analysis} {\bf 18}(1987), 208-218.

\bibitem{dieu-t}  Dieudonne, J.  {\em Treatise on
Analysis}, vol. VI, Academic Press, 1978.

\bibitem{dieudonneshort}  Dieudonne, J.  {\em Special Functions
and Linear Representations of Lie Groups}, CBMS Regional Conf. Ser.,
no. 42, AMS, 1980.

\bibitem{drmotatichy}
Drmota, M. and Tichy, R.F.
{\em Sequences, Discrepancies and Applications},
Lecture Notes in Math. 1651,
Springer, 1997.

\bibitem{dym-mckean}
Dym, H. and McKean, H.P.
{\em Fourier Series and Integrals},
Academic Press, New York, 1972.

\bibitem{gibbs-su}
Gibbs, A.L. and Su, F.E.
On choosing and bounding probability metrics, preprint.

\bibitem{gr}  Greenhalgh, A.  {\em Random walks on groups with
subgroup invariance properties}.  Ph.D. Thesis, Dept. of
Mathematics, Stanford Univ., 1987.

\bibitem{hensley-su} Hensley, D. and Su, F.E. Random walks with
badly approximable numbers, in {\em Unusual Applications of Number
Theory}, DIMACS Ser. Discrete Math. Theoret. Comput. Sci., AMS,
to appear.

\bibitem{hr}  Hewitt, E. and Ross, K.  {\em Abstract
Harmonic Analysis}, Vol. II.  Springer-Verlag, 1970.

\bibitem{j}  Jackson, D.  {\em Fourier Series and Orthogonal
Polynomials}.  The Carus Mathematical Monographs, No. 6.  MAA,
1941.

\bibitem{letac} Letac, G.  Problemes classiques de probabilite sur
un couple de Gelfand.  {\em Lecture Notes in Math.}, no. 861,
Springer-Verlag, 1981.

\bibitem{kn}
L. Kuipers and H. Neiderreiter,
{\em Uniform Distribution of Sequences},
Wiley, New York, 1974.

\bibitem{porod}  Porod, U.  {\em Time to stationarity for random walks on
compact Lie groups}.  Ph.D. thesis, Division of Math. Sciences,
Johns Hopkins University, 1993.

\bibitem{porod-ann}  Porod, U.  The cutoff phenomenon for random
reflections, {\em Ann. Probab.} {\bf 24}(1996), 74-96.

\bibitem{rosen1}  Rosenthal, J.S.  Random rotations: characters and
random walks on $SO(n)$.  {\em Ann. Probab.}, {\bf 22}(1994),
398-423.

\bibitem{suthesis} Su, F.E.  {\em Methods for Quantifying Rates of Convergence
for Random Walks on Groups.}  Ph.D.\ Thesis, Harvard University, 1995.

\bibitem{su-trans} Su, F.E.  Convergence of random walks on the circle
generated by an irrational rotation.
{\em Trans. Amer. Math. Soc.}, 350(1998), 3717-3741.

\bibitem{su-leveque} Su, F.E.  A LeVeque-type lower bound for discrepancy,
in {\em Monte Carlo and Quasi-Monte Carlo Methods 1998}, H.
Niederreiter and J. Spanier, eds., Springer-Verlag, 2000, 448-458.

\bibitem{voit-sphere} Voit, M.  A central limit theorem for
isotropic random walks on $n$-spheres for $n\rightarrow\infty$,
{\em J. Math. Anal. Appl.} {\bf 189}(1995), 215-224.

\bibitem{voit-twopt} Voit, M.  Limit theorems for compact
two-point homogeneous spaces of large dimensions,  {\em J.
Theoret. Probab.} {\bf 9}(1996), 353-370.

\bibitem{voit-ann} Voit, M.  Rate of convergence to Gaussian
measures on $n$-spheres and Jacobi Hypergroups, {\em Ann. Probab.}
{\bf 25} (1997), 457-477.


\end{thebibliography}
\end{document}